# WEAK POINCARÉ INEQUALITIES ON DOMAINS DEFINED BY BROWNIAN ROUGH PATHS[1]

By Shigeki Aida

*Osaka University*

We prove weak Poincaré inequalities on domains which are inverse images of open sets in Wiener spaces under continuous functions of Brownian rough paths. The result is applicable to Dirichlet forms on loop groups and connected open subsets of path spaces over compact Riemannian manifolds.

**1. Introduction.** Let $\mathbf{w}(t)$ be the $d$-dimensional standard Brownian motion starting at the origin. Let $\overline{\mathbf{w}}(s,t)_1 = \mathbf{w}(t) - \mathbf{w}(s)$. Also let us consider a two parameter process with values in $\mathbb{R}^d \otimes \mathbb{R}^d$ defined by a Stratonovich stochastic integral

$$\overline{\mathbf{w}}(s,t)_2 = \int_s^t (\mathbf{w}(u) - \mathbf{w}(s)) \otimes d\mathbf{w}(u), \tag{1.1}$$

where $0 \leq s \leq t \leq 1$ and $\otimes$ denotes a tensor product. Lyons [17] proved that solutions of stochastic differential equations (SDEs) are continuous functions of the Brownian rough path $\overline{\mathbf{w}}(s,t) = (\overline{\mathbf{w}}(s,t)_1, \overline{\mathbf{w}}(s,t)_2)$. We give a precise definition of the Brownian rough path in the next section; see also [18] and [15]. The discontinuity of solutions of SDEs in the uniform convergence topology of $\mathbf{w}$ causes difficulties in analysis on Wiener spaces. However, the Lyons results provide a good topology on Wiener space and may be applied to problems which have difficulties because of the discontinuity of Wiener functionals; for example, see [16]. The present paper is an attempt to apply the Lyons continuity theorem to problems in infinite-dimensional analysis and we prove weak Poincaré inequalities (WPIs) on some "connected" domain on a Wiener space defined by a continuous function of Brownian

Received January 2003; revised August 2003.

[1]Supported in part by Grant-in-Aid for Scientific Research (C) 12640173 and the Sumitomo Foundation.

*AMS 2000 subject classifications.* Primary 60H07; secondary 60H10, 58J65.

*Key words and phrases.* Weak Poincaré inequality, Brownian rough path, convexity, logarithmic Sobolev inequality.

rough paths. The WPI (actually, equivalent uniform positivity improving property of the corresponding diffusion semigroup) on a connected domain was first proved by Kusuoka [14] and led to abundant research on analysis on Wiener space and loop space. The WPI itself was introduced in [21] and the equivalence to uniform positivity improving property of the semigroup was proved therein. Aida [4] proved that WPI holds on a domain which is a "connected union" of domains on which WPI hold with respect to the natural Dirichlet forms. This proves that a WPI holds on a connected open set on a Wiener space with respect to the natural Dirichlet form because the Poincaré inequality (PI) [actually, stronger logarithmic Sobolev inequality (LSI)] holds on a ball [1, 4, 6]. Note that the inverse image of an open set by a solution of SDE is not an open set in the usual topology and the above argument is not applicable to such a set. However, by the Lyons theorem, by replacing the usual ball with a finer ball in the sense of rough path, we may apply the above argument. This is the main idea of this paper.

The structure of this paper is as follows. In Section 2, we introduce notation and state WPIs on our unit set $U_{a,\mathbf{z}}$ defined by Brownian rough path $\overline{\mathbf{w}}$, which plays the role of a ball in the usual continuous category. The domain is a nonconvex set which is defined by a quadratic Wiener functional, that is, Lévy's stochastic area in the Wiener space. However, comparing to convex case, it seems that useful criteria for the validity of LSI, PI and WPI on unbounded nonconvex domains are not known. In Section 3, we prove a general result, Lemma 3.1, which enables us to prove WPI on nonconvex domains. This is a generalization of the fact that PI is stable when taking the product of the state spaces. Using this result, we prove WPI on $U_{a,\mathbf{z}}$ by an induction on the dimension of the Wiener space. We use the validity of LSI on a convex domain (Lemma 3.4) as the first step of the induction and, next, we use Lemma 3.1 to prove general cases. In Section 4, we prove the main theorem and WPI on a domain on path spaces and a loop group. Note that key results in [1] for the proof of the weak spectral gap property, which is equivalent to the validity of WPI on loop spaces, are the existence of a good tubular neighborhood of a submanifold (which is obtained by a solution of SDE) and a good retract map. We need to show that the tubular neighborhood can be represented as a connected union of our finer domains to prove a WPI on general loop spaces over compact Riemannian manifolds. We will study the general cases and a concrete estimate on a function $\xi(\cdot)$ in WPI in separate papers.

**2. Preliminaries and notation.** Let $T_2(\mathbb{R}^d) = \mathbb{R}^d \oplus (\mathbb{R}^d \otimes \mathbb{R}^d)$. Let $C(\Delta, T_2(\mathbb{R}^d))$ be a space of continuous functions on a simplex $\Delta = \{(s,t) \in \mathbb{R}^2 | 0 \le s \le t \le 1\}$ with values in $T_2(\mathbb{R}^d)$. Let $q > 1$. For $\eta : \Delta \to \mathbb{R}$, $\|\eta\|_q$ is defined by

$$\|\eta\|_q = \sup_D \left\{ \sum_{i=0}^{n-1} |\eta(t_i, t_{i+1})|^q \right\}^{1/q}, \tag{2.1}$$

where $D = \{0 = t_0 < t_1 < \cdots < t_n = 1\}$ runs all partitions of $[0,1]$. Let $e_i = {}^t(0,\ldots,1^i,\ldots,0)$. For $\eta = (\eta(\cdot,\cdot)_1, \eta(\cdot,\cdot)_2) \in C(\Delta, T_2(\mathbb{R}^d))$, set $\eta_{1,i}(s,t) = (\eta(s,t)_1, e_i)$, $\eta_{2,k,l}(s,t) = (\eta(s,t)_2, e_k \otimes e_l)$ and define

$$\|\eta(\cdot,\cdot)_1\|_q = \max_{1 \leq i \leq d} \|\eta(\cdot,\cdot)_{1,i}\|_q, \tag{2.2}$$

$$\|\eta(\cdot,\cdot)_2\|_q = \max_{1 \leq k,l \leq d} \|\eta(\cdot,\cdot)_{2,k,l}\|_q. \tag{2.3}$$

Let $p$ be a positive number such that $2 < p < 3$ and define a $p$-variation norm $\|\eta\|_{C^p(\Delta, T_2(\mathbb{R}^d))}$ for $\eta(\cdot,\cdot) \in C(\Delta, T_2(\mathbb{R}^d))$ by

$$\|\eta(\cdot,\cdot)\|_{C^p(\Delta, T_2(\mathbb{R}^d))} = \max\{\|\eta_1\|_p, \|\eta_2\|_{p/2}\}, \tag{2.4}$$

where $C^p(\Delta, T_2(\mathbb{R}^d))$ stands for the subset of $C(\Delta, T_2(\mathbb{R}^d))$ that consists of all elements $\eta$ with $\|\eta\|_{C^p(\Delta, T_2(\mathbb{R}^d))} < \infty$. Subsequently we denote $\|\cdot\|_{C^p(\Delta, T_2(\mathbb{R}^d))}$ by $\|\cdot\|_{C^p}$ for simplicity. Also $|\cdot|, \|\cdot\|$ stand for the usual Euclidean norm unless otherwise indicated. We denote by $\mathbf{W}^d$ $(= W_1 \times \cdots \times W_d)$ the set of all continuous paths $\mathbf{w}(\cdot)$ on $[0,1]$ with values in $\mathbb{R}^d$ starting at $0$ with the Wiener measure $\mu$, where $W_i$ denotes the one-dimensional Wiener space. Let $\mathbf{H}^d$ be the Cameron–Martin subspace. For $\mathbf{h} \in \mathbf{H}^d$, let $\overline{\mathbf{h}}(s,t)_1 = \mathbf{h}(t) - \mathbf{h}(s)$ and $\overline{\mathbf{h}}(s,t)_2 = \int_s^t (\mathbf{h}(u) - \mathbf{h}(s)) \otimes d\mathbf{h}(u)$. Then $\overline{\mathbf{h}}(\cdot,\cdot) = (\overline{\mathbf{h}}(\cdot,\cdot)_1, \overline{\mathbf{h}}(\cdot,\cdot)_2) \in C^p(\Delta, T_2(\mathbb{R}^d))$ and this is called a smooth rough path. The closure of all smooth rough paths in the topology of $C^p(\Delta, T_2(\mathbb{R}^d))$ is the space of geometric rough path which we denote by $G\Omega_p(\mathbb{R}^d)$. Now we consider Brownian rough paths. We denote by $P_n\mathbf{w}$ the dyadic polygonal approximation of $\mathbf{w} \in \mathbf{W}^d$ such that

$$(P_n\mathbf{w})(t) = \mathbf{w}(t_k^n) + 2^n(\mathbf{w}(t_{k+1}^n) - \mathbf{w}(t_k^n))(t - t_k^n), \qquad t_k^n \leq t \leq t_{k+1}^n,$$

where $t_k^n = k/2^n, 0 \leq k \leq 2^n$. Note that $P_n$ is a projection operator on $\mathbf{H}^d$ such that $P_n\mathbf{H}^d \subset P_{n+1}\mathbf{H}^d$ for all $n \in \mathbb{N}$ and $\lim_{n \to \infty} P_n = I_H$ strongly. Since $P_n\mathbf{w} \in \mathbf{H}^d$, we can associate a smooth rough path $\overline{P_n\mathbf{w}} \in G\Omega_p(\mathbb{R}^d)$. The following lemma was proved in [15] and [18].

LEMMA 2.1. *For almost all $\mathbf{w}$, there exists $\overline{\mathbf{w}} \in G\Omega_p(\mathbb{R}^d)$ such that $\lim_{n \to \infty} \|\overline{P_n\mathbf{w}} - \overline{\mathbf{w}}\|_{C^p} = 0$. Moreover, the convergence is in the sense of $L^1(\mu)$.*

The limit $\overline{\mathbf{w}}(s,t) = (\overline{\mathbf{w}}(s,t)_1, \overline{\mathbf{w}}(s,t)_2)$ is called a Brownian rough path. Note that $\overline{\mathbf{w}}(s,t)_2 = \int_s^t (\mathbf{w}(u) - \mathbf{w}(s)) \otimes d\mathbf{w}(u)$ a.s. $\mathbf{w}$, where the right-hand side is the Stratonovich integral. Let $V_p(\mathbb{R}^d)$ be the closure of $\mathbf{H}^d$ with respect to the norm $\|\mathbf{h}\|_p := \|\overline{\mathbf{h}}_1\|_p$. Then $\mu(V_p(\mathbb{R}^d)) = 1$ by Lemma 2.1 and $V_p(\mathbb{R}^d)$ is a separable Banach space by part 3 of Lemma 2.2. Let

$$O_a(\mathbf{h}) = \{\mathbf{w} \in \mathbf{W}^d | \|\overline{\mathbf{w}} - \overline{\mathbf{h}}\|_{C^p} < a\}. \tag{2.5}$$

This set is a candidate of a ball-like set in the category of the continuity of Brownian rough paths. For a technical reason, we introduce a different kind of set.

Let $\mathbf{h}_1 \in C([0,1] \to \mathbb{R}^d)$ and $\mathbf{h}_2 \in \mathbf{H}^m$. We consider the Stieltjes integral

$$(2.6) \qquad C_{\mathbf{h}_1,\mathbf{h}_2}(s,t) = \int_s^t (\mathbf{h}_1(u) - \mathbf{h}_1(s)) \otimes d\mathbf{h}_2(u),$$

where $0 \le s \le t \le 1$. Of course, $C_{\mathbf{h}_1,\mathbf{h}_2}(s,t)$ is also well defined in the case where $\mathbf{h}_1 \in \mathbf{H}^m, \mathbf{h}_2 \in C([0,1] \to \mathbb{R}^d)$. Note that $\overline{\mathbf{h}}(s,t)_2 = C_{\mathbf{h},\mathbf{h}}(s,t)$. For these integrals, we use the following lemma several times.

LEMMA 2.2.  1. *For $\mathbf{h}_1 \in C([0,1] \to \mathbb{R}^d)$ and $\mathbf{h}_2 \in \mathbf{H}^m$, we have*

$$(2.7) \qquad \|C_{\mathbf{h}_1,\mathbf{h}_2}\|_{p/2} \le \|\mathbf{h}_1\|_p \|\mathbf{h}_2\|_{\mathbf{H}^m}.$$

*In the case where $\mathbf{h}_1 \in \mathbf{H}^m, \mathbf{h}_2 \in C([0,1] \to \mathbb{R}^d)$, we have*

$$(2.8) \qquad \|C_{\mathbf{h}_1,\mathbf{h}_2}\|_{p/2} \le (\|\mathbf{h}_1\|_{\mathbf{H}^m} + \|\mathbf{h}_1\|_p)\|\mathbf{h}_2\|_p.$$

2. *Let $\mathbf{w}$ and $\mathbf{z}$ be continuous paths on $\mathbb{R}^d$ and $\mathbb{R}^m$, respectively. Then $\|\overline{\mathbf{w}}_1 \otimes \overline{\mathbf{z}}_1\|_{p/2} \le \|\mathbf{w}\|_p \cdot \|\mathbf{z}\|_p.$*

3. *For any $\mathbf{h}$, $\|\mathbf{h}\|_p \le \int_0^1 |\dot{\mathbf{h}}(t)| \, dt.$*

PROOF.  1. We have

$$
\begin{aligned}
\left| \int_s^t (\mathbf{h}_1(u) - \mathbf{h}_1(s)) \otimes d\mathbf{h}_2(u) \right|^{p/2} \\
\le \left( \int_s^t |\mathbf{h}_1(u) - \mathbf{h}_1(s)| \cdot |\dot{\mathbf{h}}_2(u)| \, du \right)^{p/2} \\
\le \left( (2\varepsilon)^{-1} \int_s^t |\mathbf{h}_1(u) - \mathbf{h}_1(s)|^2 \, du + 2^{-1}\varepsilon \int_s^t |\dot{\mathbf{h}}_2(u)|^2 \, du \right)^{p/2} \\
(2.9) \quad \le 2^{-1}\varepsilon^{-p/2}(t-s)^{(p/2)-1} \int_s^t |\mathbf{h}_1(u) - \mathbf{h}_1(s)|^p \, du \\
+ 2^{-1}\varepsilon^{p/2} \left( \int_s^t |\dot{\mathbf{h}}_2(u)|^2 \, du \right)^{p/2} \\
\le 2^{-1}\varepsilon^{-p/2}(t-s)^{p/2}|\mathbf{h}_1(u_*) - \mathbf{h}_1(s)|^p + 2^{-1}\varepsilon^{p/2} \left( \int_s^t |\dot{\mathbf{h}}_2(u)|^2 \, du \right)^{p/2},
\end{aligned}
$$

where $s < u_* < t$. Let $\{t_k\}_{k=1}^n$ be a partition of $[0,1]$. Noting that

$$\sum_k \left( \int_{t_{k-1}}^{t_k} |\dot{\mathbf{h}}_2(u)|^2 \, du \right)^{p/2}$$

$$
\begin{aligned}
(2.10) \qquad &= \sum_k \left(\int_0^1 |\dot{\mathbf{h}}_2(u)|^2\, du\right)^{p/2} \left(\frac{\int_{t_{k-1}}^{t_k} |\dot{\mathbf{h}}_2(u)|^2\, du}{\int_0^1 |\dot{\mathbf{h}}_2(u)|^2\, du}\right)^{p/2} \\
&\leq \|\mathbf{h}_2\|_{\mathbf{H}^m}^p,
\end{aligned}
$$

we get

$$
(2.11) \qquad \sum_{k=1}^n \left| \int_{t_{k-1}}^{t_k} (\mathbf{h}_1(u) - \mathbf{h}_1(t_{k-1})) \otimes d\mathbf{h}_2(u) \right|^{p/2}
$$
$$
\leq 2^{-1}(\varepsilon^{-p/2}\|\mathbf{h}_1\|_p^p + \varepsilon^{p/2}\|\mathbf{h}_2\|_{\mathbf{H}^m}^p).
$$

Minimizing the function of $\varepsilon$ on the right-hand side, we get (2.7). Noting $C_{\mathbf{h}_1,\mathbf{h}_2} = -C^*_{\mathbf{h}_2,\mathbf{h}_1} + \overline{\mathbf{h}}_1 \otimes \overline{\mathbf{h}}_2$, (2.8) follows from (2.7) and part 2. Here $C^*_{\mathbf{h}_2,\mathbf{h}_1}$ stands for the transposed matrix under the natural identification between $\mathbb{R}^m \otimes \mathbb{R}^d$ and the space of $(m,d)$ matrices. We prove statement 2:

$$
(2.12) \qquad \|\overline{\mathbf{w}}(\cdot,\cdot)_1 \otimes \overline{\mathbf{z}}(\cdot,\cdot)_1\|_{p/2}^{p/2} \leq 2^{-1}(\varepsilon\|\mathbf{w}\|_p^p + \varepsilon^{-1}\|\mathbf{z}\|_p^p)
$$
$$
\leq \|\mathbf{w}\|_p^{p/2}\|\mathbf{z}\|_p^{p/2}.
$$

The proof of part 3 is similar to (2.10). □

The definition of $C_{\mathbf{h}_1,\mathbf{h}_2}$ can be extended to the Brownian path by using the Wiener integral. That is, for $\mathbf{z} \in C([0,1] \to \mathbb{R}^m)$, we can define, for almost all $\mathbf{w} \in \mathbf{W}^d$,

$$
(2.13) \qquad C_{\mathbf{w},\mathbf{z}}(s,t) = \int_s^t \mathbf{w}(s,u)_1 \otimes d\mathbf{z}(u)
$$

as the Wiener integral. Clearly, $C_{\mathbf{w},\mathbf{z}}(\cdot,\cdot) \in C(\Delta, \mathbb{R}^d \otimes \mathbb{R}^m)$. Actually, if $\mathbf{z}$ has more regularity, then so does $C_{\mathbf{w},\mathbf{z}}$. To show this, we introduce a norm which is useful in many calculations. For $\mathbf{z} \in C([0,1] \to \mathbb{R}^m)$ and $\kappa > p-1$, let

$$
(2.14) \qquad \|\mathbf{z}\|_{p,\kappa} := \left[\sum_{n=1}^\infty \left\{ n^\kappa \sum_{k=1}^{2^n} |\mathbf{z}(t_k^n) - \mathbf{z}(t_{k-1}^n)|^p \right\}\right]^{1/p}.
$$

By the Schwarz inequality, we have $\|\mathbf{z}\|_{p,\kappa} \leq C\|\mathbf{z}\|_{\mathbf{H}^m}$. Also it is easy to check that if $\|\mathbf{z}\|_{p,\kappa} < \infty$, then $\|P_N\mathbf{z}\|_{p,\kappa} < \infty$ and $\lim_{N\to\infty} \|P_N\mathbf{z} - \mathbf{z}\|_{p,\kappa} = 0$. The estimate below can be found in Lemma 2 in [15] and Proposition 4.1.1 in [18]. There exists a positive number $C$ such that

$$
(2.15) \qquad \|\mathbf{z}\|_p \leq C\|\mathbf{z}\|_{p,\kappa} \qquad \text{for all } \mathbf{z} \in C([0,1] \to \mathbb{R}^m).
$$

We denote by $V_{p,\kappa}(\mathbb{R}^m)$ the space that consists of all $\mathbf{z} \in V_p(\mathbb{R}^m)$ with $\|\mathbf{z}\|_{p,\kappa} < \infty$. By the above results, we see that $V_{p,\kappa}(\mathbb{R}^m)$ is a separable Banach space and $\mu(V_{p,\kappa}(\mathbb{R}^m)) = 1$.

We give estimates on $C_{\mathbf{w},\mathbf{z}}$.

LEMMA 2.3. *Let $\mathbf{z} \in V_{p,\kappa}(\mathbb{R}^m)$.*

1. *It holds that*

$$E[\|C_{\mathbf{w},\mathbf{z}}\|_{p/2}^{p/2}] \leq C(p,\kappa)\|\mathbf{z}\|_{p,\kappa}^{p/2}. \tag{2.16}$$

*Here $C(p,\kappa)$ is a positive constant which depends only on $p$ and $\kappa$.*

2. *It holds that*

$$\lim_{n\to\infty} E[\|C_{\mathbf{w},\mathbf{z}} - C_{P_n\mathbf{w},\mathbf{z}}\|_{p/2}] = 0. \tag{2.17}$$

PROOF. 1. We use the argument in Proposition 4.1.1 in [18] and Lemma 2 in [15]. Since

$$C_{\mathbf{w},\mathbf{z}}(s,t) = (\mathbf{w}(t) - \mathbf{w}(s)) \otimes (\mathbf{z}(t) - \mathbf{z}(s)) - \int_s^t d\mathbf{w}(u) \otimes (\mathbf{z}(u) - \mathbf{z}(s)), \tag{2.18}$$

by part 2 of Lemma 2.2, it suffices to estimate $C_{\mathbf{z},\mathbf{w}}(\cdot,\cdot)$. Note that for any partitions $D = \{s_i\}_{i=0}^N$ of $[s,t]$,

$$C_{\mathbf{z},\mathbf{w}}(s,t) = \sum_{i=1}^N C_{\mathbf{z},\mathbf{w}}(s_{i-1},s_i) + \sum_{1\leq i < j \leq N} (\mathbf{z}(s_i) - \mathbf{z}(s_{i-1})) \otimes (\mathbf{w}(s_j) - \mathbf{w}(s_{j-1}))$$

holds. Thus

$$|C_{\mathbf{z},\mathbf{w}}(s,t)|^{p/2}$$
$$\leq 2^{(p/2)-1} \left|\sum_{i=1}^N C_{\mathbf{z},\mathbf{w}}(s_{i-1},s_i)\right|^{p/2}$$
$$+ 2^{(p/2)-1} \left|\sum_{1\leq i < j \leq N} (\mathbf{z}(s_i) - \mathbf{z}(s_{i-1})) \otimes (\mathbf{w}(s_j) - \mathbf{w}(s_{j-1}))\right|^{p/2} \tag{2.19}$$
$$\leq 2^{(p/2)-1} \left\{\left|\sum_{i=1}^N C_{\mathbf{z},\mathbf{w}}(s_{i-1},s_i)\right|^{p/2} + \varepsilon^{-1}\left(\sum_{i=1}^N |\mathbf{z}(s_i) - \mathbf{z}(s_{i-1})|\right)^p\right.$$
$$\left.+ \varepsilon\left(\sum_{i=1}^N |\mathbf{w}(s_i) - \mathbf{w}(s_{i-1})|\right)^p\right\},$$

where $\varepsilon$ is a positive number. Since $\kappa > p - 1 > (p/2) - 1$, by the same argument as in the proof of Proposition 4.1.1 in [18], we have

$$\|C_{\mathbf{z},\mathbf{w}}\|_{p/2}^{p/2} \leq C(p,\kappa) \sum_{n=1}^\infty n^\kappa \sum_{k=1}^{2^n} (|C_{\mathbf{z},\mathbf{w}}(t_{k-1}^n,t_k^n)|^{p/2} + \varepsilon^{-1}|\mathbf{z}(t_k^n) - \mathbf{z}(t_{k-1}^n)|^p$$
$$+ \varepsilon|\mathbf{w}(t_k^n) - \mathbf{w}(t_{k-1}^n)|^p), \tag{2.20}$$

where $C(p,\kappa)$ is a constant which depends only on $p$ and $\kappa$, although constants may change line by line in the calculation below. We estimate the expectation of $C_{\mathbf{z},\mathbf{w}}(\cdot,\cdot)$,

$$E[|C_{\mathbf{z},\mathbf{w}}(t_{k-1}^n, t_k^n)|^{p/2}] \leq C\left\{\int_{t_{k-1}^n}^{t_k^n} |\mathbf{z}_u - \mathbf{z}_{t_{k-1}^n}|^2 \, du\right\}^{p/4}$$

(2.21)

$$= C|\mathbf{z}_{u_k} - \mathbf{z}_{t_{k-1}^n}|^{p/2}\left(\frac{1}{2^n}\right)^{p/4},$$

where $t_{k-1}^n \leq u_k \leq t_k^n$. Thus we get

(2.22)
$$E[|C_{\mathbf{z},\mathbf{w}}(t_{k-1}^n, t_k^n)|^{p/2}]$$
$$\leq \frac{C}{2}\left\{\varepsilon^{-1}(|\mathbf{z}_{t_k^n} - \mathbf{z}_{u_k}|^p + |\mathbf{z}_{u_k} - \mathbf{z}_{t_{k-1}^n}|^p)\left(\frac{1}{2^n}\right)^{\delta/2} + \varepsilon\left(\frac{1}{2^n}\right)^{(p-\delta)/2}\right\},$$

where $0 < \delta < p-2$. By using this estimate and taking the expectation of both sides in (2.20), we have

(2.23)
$$E[\|C_{\mathbf{z},\mathbf{w}}\|_{p/2}^{p/2}] \leq C(p,\kappa)\left[\sum_{n=1}^{\infty} n^\kappa \left\{\varepsilon^{-1}\|\mathbf{z}\|_p^p\left(\frac{1}{2^n}\right)^{\delta/2} + \varepsilon\left(\frac{1}{2^n}\right)^{(p-\delta-2)/2}\right\}\right.$$
$$\left. + \varepsilon^{-1}\|\mathbf{z}\|_{p,\kappa}^p + \varepsilon\sum_{n=1}^{\infty} n^\kappa\left(\frac{1}{2^n}\right)^{(p/2)-1}\right].$$

Therefore, by choosing $\varepsilon$ to minimize the right-hand side, we get the desired estimate.

2. Let $P_n^\perp \mathbf{w} = \mathbf{w} - P_n\mathbf{w}$. Then it is easy to see that for any $q > 1$, $E[|P_n^\perp \mathbf{w}(t) - P_n^\perp \mathbf{w}(s)|^q] \leq C|t-s|^{q/2}$, where $C$ is a positive number independent of $n$, $t$ and $s$. Also $E[|C_{P_n^\perp \mathbf{w},\mathbf{z}}(s,t)|^2] \leq C \cdot E[|C_{\mathbf{w},\mathbf{z}}(s,t)|^2]$ for any $n$, $E[|C_{P_n^\perp \mathbf{w},\mathbf{z}}(s,t)|^{p/2}] \leq C \cdot E[|C_{P_n^\perp \mathbf{w},\mathbf{z}}(s,t)|^2]^{p/4}$, $P_n^\perp \mathbf{w}(t_m^k) = 0$ for all $m \leq n$ and $\lim_{n\to\infty} E[|C_{P_n^\perp \mathbf{w},\mathbf{z}}(s,t)|^2] = 0$ for all $t,s$. Hence by an argument similar to statement 1, we can complete the proof of part 2. $\square$

Let us introduce a subset of $C^q(\Delta, \mathbb{R}^d \otimes \mathbb{R}^m)$ $(1 < q < \frac{3}{2})$, $V_q(\Delta, \mathbb{R}^d \otimes \mathbb{R}^m)$, which is the closure of the following linear subspace in the $q$-variation norm:

(2.24)
$$\left\{\eta \in C^q(\Delta, \mathbb{R}^d \otimes \mathbb{R}^m) \,\Big|\, \eta = \sum_{i=1}^n C_{\varphi_i, \phi_i},\right.$$
$$\left. \text{where } \varphi_i \in \mathbf{H}^d, \phi_i \in \mathbf{H}^m \text{ and } n \in \mathbb{N}\right\}.$$

Lemma 2.2 implies that $V_q(\Delta, \mathbb{R}^d \otimes \mathbb{R}^m)$ is a separable Banach space. By Lemmas 2.2 and 2.3, we can find a version of $C_{\mathbf{w},\mathbf{z}}$ with values in $V_{p/2}(\Delta, \mathbb{R}^d \otimes$ $\mathbb{R}^m$) and an $\mathbf{H}^d$-invariant subset of $\mathbf{W}^d$ such that $\mathbf{X}_1^d + \mathbf{H}^d = \mathbf{X}_1^d$ with $\mu(\mathbf{W}^d \setminus \mathbf{X}_1^d) = 0$ and for all $\mathbf{w} \in \mathbf{X}_1^d$ and $\mathbf{h} \in \mathbf{H}^d$,

$$(2.25) \quad C_{\mathbf{w}+\mathbf{h},\mathbf{z}}(s,t) = C_{\mathbf{w},\mathbf{z}}(s,t) + C_{\mathbf{h},\mathbf{z}}(s,t) \qquad \text{for all } 0 \leq s \leq t \leq 1.$$

We note that by Lemma 2.2, $C_{\mathbf{w},\mathbf{z}}$ is a $V_{p/2}(\Delta, \mathbb{R}^d \otimes \mathbb{R}^m)$-valued $\mathbf{H}^d$-continuous function. Note that $C_{\mathbf{z},\mathbf{w}} := \bar{\mathbf{z}}_1 \otimes \overline{\mathbf{w}}_1 - C_{\mathbf{w},\mathbf{z}}^*$ is a version of Wiener integral $\int_s^t (\mathbf{z}(u) - \mathbf{z}(s)) \otimes d\mathbf{w}(u)$.

We consider another $\mathbf{H}^d$-invariant subset. Let $\mathbf{X}_2^d$ be the set of all $\mathbf{w}$ which converge as stated in Lemma 2.1. Then $\mathbf{X}_2^d$ is $\mathbf{H}^d$ invariant. To show this note that for any $\varphi, \phi \in \mathbf{H}^d$,

$$\begin{aligned}
&\bar{\varphi}(s,t)_2 - \bar{\phi}(s,t)_2 \\
&= \int_s^t [(\varphi - \phi)(u) - (\varphi - \phi)(s)] \otimes d(\varphi(u) - \phi(u)) \\
&\quad + \int_s^t [(\varphi - \phi)(u) - (\varphi - \phi)(s)] \otimes d\phi(u) \\
&\quad - \int_s^t d\phi(u) \otimes [(\varphi - \phi)(u) - (\varphi - \phi)(s)] \\
&\quad + (\phi(t) - \phi(s)) \otimes [(\varphi - \phi)(t) - (\varphi - \phi)(s)].
\end{aligned} \quad (2.26)$$

Putting $\varphi = P_n \mathbf{w}$, $\phi = P_n \mathbf{h}$ and combining Lemma 2.2, $\overline{P_n(\mathbf{w} - \mathbf{h})}$ converges in $C^p$ for any $\mathbf{w} \in \mathbf{X}_2^d$ and $\mathbf{h} \in \mathbf{H}^d$. In this sense, $\overline{\mathbf{w} - \mathbf{h}}$ is well defined and (2.26) still holds for $\varphi = \mathbf{w} \in \mathbf{X}_2^d$. Also $\mathbf{w}(\in \mathbf{X}_2^d) \to \overline{\mathbf{w}}(\in C^{p/2}(\Delta, \mathbb{R}^d \otimes \mathbb{R}^d))$ is an $\mathbf{H}^d$-continuous function by this equation and Lemma 2.2.

We need results analogous to Lemma 2.1 for our purposes.

LEMMA 2.4. *For almost all $\mathbf{w}$ and in the sense of $L^1$ as $N \to \infty$, $\|\overline{(\mathbf{w} - P_N \mathbf{w})}_2\|_{p/2}$ and $\|C_{\mathbf{w} - P_N \mathbf{w}, P_N \mathbf{w}}\|_{p/2}$ converge to $0$.*

PROOF. We use (2.26) in the case where $\varphi = \mathbf{w}$ and $\phi = P_N \mathbf{w}$. Note that $E[|P_N \mathbf{w}(t) - P_N \mathbf{w}(s)|^p] \leq C|t - s|^{p/2}$. Here $C$ is a positive number independent of $N$. Hence, $\sup_N E[\|P_N \mathbf{w}\|_{p,\kappa}^p] < \infty$. Also by Lemma 2.1 and the independent property of $P_N \mathbf{w}$ and $P_N^\perp \mathbf{w} = \mathbf{w} - P_N \mathbf{w}$, the fourth term on the right-hand side of (2.26) can be estimated by Lemma 2.2.2 and the Itô–Nisio theorem [12]. So it suffices to prove that $\|C_{P_N^\perp \mathbf{w}, P_N \mathbf{w}}\|_{p/2}$ converges to $0$. Using the independence of $P_N \mathbf{w}$ and $P_N^\perp \mathbf{w}$ and the Gaussian property, we have

$$(2.27) \; E[|C_{P_N^\perp \mathbf{w}, P_N \mathbf{w}}(t_{k-1}^n, t_k^n)|^{p/2}] \leq CE[|P_N^\perp \mathbf{w}(u_k) - P_N^\perp \mathbf{w}(t_{k-1}^n)|^{p/2}]\left(\frac{1}{2^n}\right)^{p/4},$$

where $t_{k-1}^n \leq u_k \leq t_k^n$. By this, we have

$$(2.28) \quad E\left[\sum_{n=1}^{\infty} n^\kappa \left(\sum_{k=1}^{2^n} |C_{P_N^\perp \mathbf{w}, P_N \mathbf{w}}(t_{k-1}^n, t_k^n)|^{p/2}\right)\right] \leq CE[\|P_N^\perp \mathbf{w}\|_p^{p/2}]$$
$$\leq C_{p,\kappa} E[\|P_N^\perp \mathbf{w}\|_{p,\kappa}^{p/2}].$$

Noting that $P_N^\perp \mathbf{w}(t_m^k) = 0$ for all $1 \leq m \leq N$ and $0 \leq k \leq 2^m$ noting that $E[|P_N^\perp \mathbf{w}(t) - P_N^\perp \mathbf{w}(s)|^p] \leq C|t-s|^{p/2}$ ($C$ is independent of $N$), and using Lemma 2 in [15], we get $E[\|P_N^\perp \mathbf{w}\|_{p,\kappa}^p] \leq 2^{-rN}$, where $r$ is a small positive number. Hence by the same argument as in the proof of Lemma 2.3.1 and the Borel–Cantelli lemma, $\overline{(\mathbf{w} - P_N \mathbf{w})}_2$ converges to 0 in $L^1$ and for almost all $\mathbf{w}$.

$\square$

By (2.26) and Lemma 2.2, it is easy to see that the set $\mathbf{X}_3^d$ of all $\mathbf{w}$ for which the convergences in Lemma 2.4 are valid is $\mathbf{H}^d$ invariant. Let $\mathbf{X}^d = \bigcap_{i=1}^3 \mathbf{X}_i^d$. Clearly, $\mathbf{X}^d$ is also an $\mathbf{H}^d$-invariant subset with $\mu(\mathbf{X}^d) = 1$.

Now we define our unit set for $a > 0$ and $\mathbf{z} \in V_{p,\kappa}(\mathbb{R}^d)$:

$$(2.29) \quad U_{a,\mathbf{z}} = \{\mathbf{w} \in \mathbf{X}^d | \|\overline{\mathbf{w}}(\cdot, \cdot)\|_{C^p} < a, \|C_{\mathbf{w},\mathbf{z}}\|_{p/2} < a, \|C_{\mathbf{z},\mathbf{w}}\|_{p/2} < a\}.$$

Let $d\mu_{a,\mathbf{z}} = d\mu|_{U_{a,\mathbf{z}}} / \mu(U_{a,\mathbf{z}})$. We prove $\mu(U_{a,\mathbf{z}}) > 0$ in Lemma 2.6. Since $U_{a,\mathbf{z}}$ is an $\mathbf{H}^d$-open set with positive measure, we can define a Dirichlet form $[\mathcal{E}_{a,\mathbf{z}}, \mathrm{D}(\mathcal{E}_{a,\mathbf{z}})]$ on $L^2(U_{a,\mathbf{z}}, d\mu_{a,\mathbf{z}})$. It is the smallest closed extension of

$$(2.30) \quad \mathcal{E}_{a,\mathbf{z}}(f,f) = \int_{U_{a,\mathbf{z}}} |Df(\mathbf{w})|_H^2 \, d\mu_{a,\mathbf{z}} \quad \text{for all } f \in \mathfrak{FC}_b^\infty|_{U_{a,\mathbf{z}}}.$$

Here $\mathfrak{FC}_b^\infty$ is the set of smooth cylindrical functions with bounded derivatives. More generally, we can define a closeable Dirichlet form on the $\mathbf{H}^d$-open set domain $U$ such that $\mathcal{E}_U(f,f) = \int_U |Df(\mathbf{w})|_H^2 \, d\mu_U(\mathbf{w})$, where $f \in \mathfrak{FC}_b^\infty|_U$ and $d\mu_U(\cdot) = d\mu(\cdot)/\mu(U)$. Refer to [3] and the references therein for the definition. The Dirichlet form is independent of the choice of the $\mathbf{H}^d$-continuous version of the defining functions of the domain. In this paper, we consider the smallest closed extension only. We prove WPI for $\mathcal{E}_{a,\mathbf{z}}$ in the next section; that is, we prove the following lemma.

LEMMA 2.5. *There exists a nonnegative jointly measurable function $\xi(\delta, a, \mathbf{z})$ $[(\delta, a, \mathbf{z}) \in (0, \infty) \times (0, \infty) \times V_{p,\kappa}(\mathbb{R}^d)]$ which is a nonincreasing function of $\delta$ such that for any $f \in \mathrm{D}(\mathcal{E}_{a,\mathbf{z}})$,*

$$(2.31) \quad \int_{U_{a,\mathbf{z}}} (f(\mathbf{w}) - \langle f \rangle_{\mu_{a,\mathbf{z}}})^2 \, d\mu_{a,\mathbf{z}}(\mathbf{w}) \leq \xi(\delta, a, \mathbf{z}) \mathcal{E}_{a,\mathbf{z}}(f,f) + \delta \|f\|_\infty^2,$$

*where $\langle f \rangle_{\mu_{a,\mathbf{z}}}$ stands for the expectation with respect to $\mu_{a,\mathbf{z}}$.*

The joint measurability is almost obvious by the definition of $\xi(\delta, a, \mathbf{z})$. The problem is to prove the boundedness. Finally, we prove the positivity of the measure of $U_{a,\mathbf{z}}$ for all $a > 0$ and $\mathbf{z}$ although it seems to be almost obvious. Note that we will use the notations in the proof of Lemma 2.6 in the argument below.

LEMMA 2.6. *For all $a > 0$ and $\mathbf{z} \in V_{p,\kappa}(\mathbb{R}^m)$, $\mu(U_{a,\mathbf{z}}) > 0$.*

PROOF. We prove this by induction on the dimension. We consider the case where $d = 2$. Let

$$(2.32) \quad U_{a,\mathbf{z},1} = \{w_1 \in W_1 | \|w_1\|_p < a, \|C_{w_1,\mathbf{z}}\|_{p/2} < a, \|C_{\mathbf{z},w_1}\|_{p/2} < a\},$$

$$(2.33) \quad \begin{aligned} W_{a,(w_1,\mathbf{z}),2} &= \{w_2 \in W_2 | \|C_{w_2,(w_1,\mathbf{z})}\|_{p/2} < a, \\ &\|C_{(w_1,\mathbf{z}),w_2}\|_{p/2} < a, \|w_2\|_p < a\}. \end{aligned}$$

Note that a measure $\nu$ on a Banach space $B$ is called a Gaussian measure with mean 0 if all random variables $\varphi(x)$ are one-dimensional Gaussian random variables with mean 0, where $\varphi \in B^*$ and $x \in B$. Although we do not determine the dual spaces $V_{p/2}(\Delta, \mathbb{R}^m)$ of $V_p(\mathbb{R})$, we can conclude that the law of $(w_1, C_{w_1,\mathbf{z}}, C_{\mathbf{z},w_1})$ defines a Gaussian measure on $V_p(\mathbb{R}) \times V_{p/2}(\Delta, \mathbb{R}^m) \times V_{p/2}(\Delta, \mathbb{R}^m)$ with mean 0. Let us explain it. For $C_{w_1,\mathbf{z}}$, there exist $\sigma(w_1)$-measurable independent Gaussian random variables with mean 0, $\{\xi_n\}_n$ and $\{C_n(\cdot, \cdot)\} \subset V_{p/2}(\Delta, \mathbb{R}^m)$, and a subsequence $\{n(k)\}_{k=1}^\infty \subset \mathbb{N}$ such that $C_{w_1,\mathbf{z}} = \lim_{k \to \infty} \sum_{l=1}^{n(k)} C_l \xi_l$ for almost all $w_1$ in the topology of $V_{p/2}(\Delta, \mathbb{R}^m)$. Clearly the law of $\sum_{l=1}^{n(k)} C_l \xi_l$ is a Gaussian measure with mean 0 on $V_{p/2}(\Delta, \mathbb{R}^m)$. These statements imply that the law of $C_{w_1,\mathbf{z}}$ is also a Gaussian measure with mean 0. The proof of the Gaussian property of other random variables is similar. Generally, any neighborhood of a 0 vector has a positive measure for any Gaussian measure with mean 0 on a separable Banach space; see Theorem 3.6.1 in [5]. Thus, $\mu(U_{a,\mathbf{z},1}) > 0$. The set (2.33) is defined for almost all $w_2$ for each $w_1 \in V_{p,\kappa}(\mathbb{R})$. Also the law of $(C_{(w_1,\mathbf{z}),w_2}, C_{w_2,(w_1,\mathbf{z})}, w_2)$ is a Gaussian measure on $V_{p/2}(\Delta, \mathbb{R}^{m+1}) \times V_{p/2}(\Delta, \mathbb{R}^{m+1}) \times V_p(\mathbb{R})$ with mean vector 0 for each fixed $w_1$, so $\mu(W_{a,(w_1,\mathbf{z}),2}) > 0$ for almost all $w_1 \in V_{p,\kappa}(\mathbb{R})$. Since $U_{a,\mathbf{z}}$ coincides with

$$\{(w_1, w_2) \in \mathbf{W}^2 | w_1 \in U_{a,\mathbf{z},1}, w_2 \in W_{a,(w_1,\mathbf{z}),2}\},$$

except a null set, by the Fubini theorem, we have $\mu(U_{a,\mathbf{z}}) > 0$. Next we prove the $(d+1)$-dimensional case by using the $d$-dimensional case. We denote $\mathbf{w} = (\mathbf{w}', w_{d+1}) \in \mathbf{W}^d \times W_{d+1}$. For a given $\mathbf{z}$, we consider a domain $U_{a,\mathbf{z}}$ in $\mathbf{W}^{d+1}$. Let

$$(2.34) \quad U_{a,\mathbf{z},d} = \{\mathbf{w}' \in \mathbf{W}^d | \|\overline{\mathbf{w}}'\|_{C^p} < a, \|C_{\mathbf{w}',\mathbf{z}}\|_{p/2} < a, \|C_{\mathbf{z},\mathbf{w}'}\|_{p/2} < a\},$$

$$W_{a,(\mathbf{w}',\mathbf{z}),d+1} = \{w_{d+1} \in W_{d+1} | \|C_{(\mathbf{w}',\mathbf{z}),w_{d+1}}\|_{p/2} < a,$$
(2.35)
$$\|C_{w_{d+1},(\mathbf{w}',\mathbf{z})}\|_{p/2} < a, \|w_{d+1}\|_p < a\}.$$

Then, for almost all $\mathbf{w}$,

(2.36) $\quad U_{a,\mathbf{z}} = \{\mathbf{w} \in \mathbf{W}^{d+1} | \mathbf{w}' \in U_{a,\mathbf{z},d}, w_{d+1} \in W_{a,(\mathbf{w}',\mathbf{z}),d+1}\}.$

By the same reasoning as in the case where $d=1$, $\mu(W_{a,(\mathbf{w}',\mathbf{z}),d+1}) > 0$ for almost all $\mathbf{w}' \in V_{p,\kappa}(\mathbb{R}^d)$. Therefore, we complete the proof by the Fubini theorem again. □

**3. WPI on $U_{a,\mathbf{z}}$.** We begin by proving a lemma in general settings which is used to prove Lemma 2.5. Let $(Y_i, \mathcal{F}_i, m_i)$ $(i=1,2)$ be complete probability spaces and Dirichlet forms $[\mathcal{E}_i, D(\mathcal{E}_i)]$ on them. Let $\Gamma_i(\cdot,\cdot)$ be the square field operator of $\mathcal{E}_i$. Let us consider a completed product probability space $Y_1 \times Y_2$ and let $U$ be a measurable subset of $Y_1 \times Y_2$. For $x \in Y_1$ and $y \in Y_2$, define $U_x = \{y \in Y_2 | (x,y) \in U\}$ and $U^y = \{x \in Y_1 | (x,y) \in U\}$. Let $U_1 = \{x \in Y_1 | m_2(U_x) > 0\}$ and $U_2 = \{y \in Y_2 | m_1(U^y) > 0\}$. Let $m_x$ ($m^y$) be the normalized probability measure on the section $U_x$ ($U^y$). We consider a pre-Dirichlet form on a section $U_x$ for almost all $x \in U_1$:

(3.1) $$\mathcal{E}_{2,x}(f,f) := \int_{U_x} \Gamma_2(f,f)(y) \, dm_x(y).$$

We define $\mathcal{E}_{1,y}$ on $L^2(U^y, dm^y(x))$ in the same way. Let $f(x,y)$ be a measurable function on $Y_1 \times Y_2$ such that $f(x,\cdot) \in D(\mathcal{E}_2)$ and $f(\cdot,y) \in D(\mathcal{E}_1)$ for fixed $x,y$. For such $f$, set

(3.2) $\quad \Gamma(f,f)(x,y) = \Gamma_1(f(\cdot,y), f(\cdot,y))(x) + \Gamma_2(f(x,\cdot), f(x,\cdot))(y)$

and define

(3.3) $$\mathcal{E}_U(f,f) := \int_U \Gamma(f,f)(x,y) \, dm_U(x,y),$$

where $dm_U$ denotes the normalized probability measure of the restriction of the product measure $dm := dm_1 \otimes dm_2$ to $U$. We denote by $\mathcal{D}_U$ the set of all functions $f$ with $\mathcal{E}_U(f,f) < \infty$.

LEMMA 3.1. *Assume that the following statements hold:*

A1. *For almost all $x \in U_1$ an $y \in U_2$, there exist jointly measurable functions $\xi_2(x,\delta)$ and $\xi_1(y,\delta)$ which are nonincreasing functions of $\delta > 0$ such that WPI holds on almost all sections:*

(3.4) $$\int_{U_x} (f(y) - \langle f \rangle_{m_x})^2 \, dm_x(y) \le \xi_2(x,\delta)\mathcal{E}_{2,x}(f,f) + \delta\|f\|_\infty^2,$$

(3.5) $$\int_{U^y} (f(x) - \langle f \rangle_{m^y})^2 \, dm^y(x) \le \xi_1(y,\delta)\mathcal{E}_{1,y}(f,f) + \delta\|f\|_\infty^2.$$

A2. *For any $\varepsilon > 0$, there exist a measurable subset $U_{1,\varepsilon} \subset U_1$ and $\delta(\varepsilon) > 0$ such that*

$$m_1(U_1 \setminus U_{1,\varepsilon}) \leq \varepsilon, \tag{3.6}$$

$$m_2(U_x \cap U_{x'}) \geq \delta(\varepsilon) \quad \text{for any } x, x' \in U_{1,\varepsilon}. \tag{3.7}$$

*Then WPI holds for the pre-Dirichlet form $\mathcal{E}_U$ on the domain $\mathcal{D}_U$.*

PROOF. First we prove that for any $\varepsilon > 0$, there exists $\widehat{U}_{i,\varepsilon} \subset U_i$ with $m_i(U_i \setminus \widehat{U}_{i,\varepsilon}) \leq \varepsilon$ such that

$$\xi(\delta, \varepsilon) = \sup\{\xi_1(y, \delta), \xi_2(x, \delta) | x \in \widehat{U}_{1,\varepsilon}, y \in \widehat{U}_{2,\varepsilon}\} < \infty \quad \text{for all } \delta > 0. \tag{3.8}$$

Let $n \in \mathbb{N}$. Take $U_{n,i} \subset U_i$ such that $m_i(U_i \setminus U_{n,i}) \leq 1/n^2$ and

$$\sup\left\{\xi_1\left(y, \frac{1}{n}\right), \xi_2\left(x, \frac{1}{n}\right) \bigg| x \in U_{n,1}, y \in U_{n,2}\right\} < \infty.$$

Then, for sufficiently large $N$, it suffices to set $U_{i,\varepsilon} = \bigcap_{n=N}^{\infty} U_{n,i}$ for our purpose. Note that

$$\iint_{U \times U} (f(x,y) - f(x',y'))^2 \, dm(x,y) \, dm(x',y')$$
$$\leq \iint_{\widetilde{U}_\varepsilon} (f(x,y) - f(x',y'))^2 \, dm(x,y) \, dm(x',y') + 8\varepsilon \|f\|_\infty^2, \tag{3.9}$$

where $\widetilde{U}_\varepsilon = \{(x,y) \in U, (x',y') \in U | x \in U_{1,\varepsilon}, x' \in U_{1,\varepsilon}\}$. Let $z \in U_x \cap U_{x'}$ for $x, x' \in U_{1,\varepsilon}$. Noting that

$$I((x,y),(x',y'))$$
$$:= (f(x,y) - f(x',y'))^2$$
$$\leq 3\{(f(x,y) - f(x,z))^2 \tag{3.10}$$
$$\quad + (f(x,z) - f(x',z))^2 + (f(x',z) - f(x',y'))^2\}$$

and by assumption A2,

$$I((x,y),(x',y')) \leq \frac{3}{\delta(\varepsilon)} \int_{U_x \cap U_{x'}} (f(x,y) - f(x,z))^2 \, dm_2(z)$$
$$+ \frac{3}{\delta(\varepsilon)} \int_{U_x \cap U_{x'}} (f(x,z) - f(x',z))^2 \, dm_2(z) \tag{3.11}$$
$$+ \frac{3}{\delta(\varepsilon)} \int_{U_x \cap U_{x'}} (f(x',z) - f(x',y'))^2 \, dm_2(z)$$
$$:= I_1 + I_2 + I_3.$$

Let $\varepsilon'$ be a positive number. By using A1 and the property of $\widehat{U}_{1,\varepsilon'}$,

$$\iint_{\widetilde{U}_\varepsilon} I_1 \, dm(x,y) \, dm(x',y')$$

$$\leq \frac{3}{\delta(\varepsilon)} \int_{x \in U_{1,\varepsilon} \cap \widehat{U}_{1,\varepsilon'}} \left( \iint_{z,y \in U_x} (f(x,y) - f(x,z))^2 \, dm_2(z) \, dm_2(y) \right) dm_1(x)$$

(3.12)

$$+ \frac{12\varepsilon'}{\delta(\varepsilon)} \|f\|_\infty^2$$

$$\leq \frac{6\xi(\delta,\varepsilon')}{\delta(\varepsilon)} \int_{x \in U_{1,\varepsilon} \cap \widehat{U}_{1,\varepsilon'}} \left( \int_{y \in U_x} \Gamma_2(f(x,\cdot))(y) \, dm_2(y) \right) dm_1(x)$$

$$+ \left( \frac{12\varepsilon'}{\delta(\varepsilon)} + \frac{6\delta m(U)}{\delta(\varepsilon)} \right) \|f\|_\infty^2$$

$$= \frac{6\xi(\delta,\varepsilon')}{\delta(\varepsilon)} \iint_U \Gamma_2(f(x,\cdot))(y) \, dm_1(x) \, dm_2(y)$$

$$+ \left( \frac{12\varepsilon'}{\delta(\varepsilon)} + \frac{6\delta m(U)}{\delta(\varepsilon)} \right) \|f\|_\infty^2.$$

Next, we estimate the integral of $I_2$:

$$\iint_{\widetilde{U}_\varepsilon} I_2 \, dm(x,y) \, dm(x',y')$$

$$\leq \frac{3}{\delta(\varepsilon)} \int_{z \in \widehat{U}_{2,\varepsilon'}} \left( \iint_{x,x' \in U_z} (f(x,z) - f(x',z))^2 \, dm_1(x) \, dm_1(x') \right) dm_2(z)$$

(3.13) $\quad + \frac{12\varepsilon'}{\delta(\varepsilon)} \|f\|_\infty^2$

$$\leq \frac{6\xi(\delta,\varepsilon')}{\delta(\varepsilon)} \iint_U \Gamma_1(f(\cdot,y))(x) \, dm_1(x) \, dm_2(y)$$

$$+ \left( \frac{12\varepsilon'}{\delta(\varepsilon)} + \frac{6\delta m(U)}{\delta(\varepsilon)} \right) \|f\|_\infty^2.$$

The integral of $I_3$ can be estimated in the same way as $I_1$. Consequently we have

$$\iint_{U \times U} (f(x,y) - f(x',y'))^2 \, dm(x,y) \, dm(x',y')$$

(3.14) $\qquad \leq \frac{18\xi(\delta,\varepsilon')}{\delta(\varepsilon)} \iint_U \Gamma(f,f)(x,y) \, dm(x,y)$

$$+ \left(8\varepsilon + \frac{36\varepsilon'}{\delta(\varepsilon)} + \frac{18\delta m(U)}{\delta(\varepsilon)}\right) \|f\|_\infty^2.$$

This completes the proof. □

REMARK 3.2. Assume PI hold in A1, and assume the coefficients $\xi_2(x)$ and $\xi_1(y)$ can be taken independently of $x$ and $y$. Further assume that there exists $U_{1,0} \subset U_1$ such that $m_1(U_1 \setminus U_{1,0}) = 0$ and $\inf_{x,x' \in U_{1,0}} m_2(U_x \cap U_{x'}) > 0$. Then we see that PI holds on $U$ by the above method. The following domain is such an example. Let $Y_1 = W_1$ and $Y_2 = W_2$, that is, one-dimensional Wiener spaces. For positive numbers $a$ and $b$ with $a < b^2$, let

(3.15) $\quad U_{a,b} = \{(w_1, w_2) | \|w_1\|_p \|w_2\|_p < a, \|w_1\|_p < b, \|w_2\|_p < b\} \subset Y_1 \times Y_2.$

By Lemma 3.4, this example satisfies the above assumptions.

To apply Lemma 3.1 to our problem, we use the following lemma.

LEMMA 3.3. *For* $\mathbf{w} \in \mathbf{W}^{d+1}$, *we denote* $\mathbf{w} = (\mathbf{w}', w_{d+1}) \in \mathbf{W}^d \times W_{d+1}$. *Let* $W_{a,(\mathbf{w}',\mathbf{z}),d+1}$ *be the set given in* (2.35). *Let* $r$ *be a positive number less than* $\frac{1}{3}$. *The following estimates hold.*

1. *For any* $\varepsilon > 0$, *we have*

$$\mu(\|w_{d+1}\|_{p,\kappa} < \varepsilon, \|C_{w_{d+1},\mathbf{z}}\|_{p/2} < \varepsilon, \|C_{\mathbf{z},w_{d+1}}\|_{p/2} < \varepsilon) > 0.$$

2. *For* $0 < \varepsilon \leq C(p,\kappa)^{-2/p} a(\alpha_{a,a,\mathbf{z}} r)^{4/p}$, *it holds that*

(3.16) $\quad \mu_{\varepsilon,a,\mathbf{z}}(\max\{\|C_{\mathbf{w}',w_{d+1}}\|_{p/2}, \|C_{w_{d+1},\mathbf{w}'}\|_{p/2}\} \geq a) \leq (\alpha_{a,a,\mathbf{z}} r)^2,$

*where* $\mu_{\varepsilon,a,\mathbf{z}}$ *denotes the conditional probability measure*

(3.17) $\quad \mu_{\varepsilon,a,\mathbf{z}}(\cdot) = \mu(\cdot | \|w_{d+1}\|_{p,\kappa} < \varepsilon, \|C_{w_{d+1},\mathbf{z}}\|_{p/2} < a, \|C_{\mathbf{z},w_{d+1}}\|_{p/2} < a),$

$\alpha_{a,a,\mathbf{z}} = \mu(\|\overline{\mathbf{w}}'\|_{C^p} < a, \|C_{\mathbf{w}',\mathbf{z}}\|_{p/2} < a, \|C_{\mathbf{z},\mathbf{w}'}\|_{p/2} < a)$ *and* $C(p,\kappa)$ *is a constant which depends on* $p$ *and* $\kappa$ *only. We take* $\varepsilon$ *in the above interval in parts 3 and 4.*

3. *Define*

(3.18) $\quad V_{p,r,a,\mathbf{z},\varepsilon} = \{\mathbf{w}' \in \mathbf{W}^d | \mu_{\varepsilon,a,\mathbf{z}}(W_{a,(\mathbf{w}',\mathbf{z}),d+1}) \geq 1 - r\alpha_{a,a,\mathbf{z}}\}.$

*Then it holds that*

(3.19) $\quad \mu(V_{p,r,a,\mathbf{z},\varepsilon} | \|\overline{\mathbf{w}}'\|_{C^p} < a, \|C_{\mathbf{w}',\mathbf{z}}\|_{p/2} < a, \|C_{\mathbf{z},\mathbf{w}'}\|_{p/2} < a) \geq 1 - r.$

4. *For any* $\mathbf{w}'_1, \mathbf{w}'_2 \in V_{p,r,a,\mathbf{z},\varepsilon}$,

(3.20) $\quad \mu(W_{a,(\mathbf{w}'_1,\mathbf{z}),d+1} \cap W_{a,(\mathbf{w}'_2,\mathbf{z}),d+1}) \geq \frac{1}{3}\tilde{\alpha}_{\varepsilon,a,\mathbf{z}},$

*where* $\tilde{\alpha}_{\varepsilon,a,\mathbf{z}} = \mu(\|w_{d+1}\|_{p,\kappa} < \varepsilon, \|C_{w_{d+1},\mathbf{z}}\|_{p/2} < a, \|C_{\mathbf{z},w_{d+1}}\|_{p/2} < a).$

PROOF. 1. The law of $(w_{d+1}, C_{w_{d+1},\mathbf{z}}, C_{\mathbf{z},w_{d+1}}) \in V_{p,\kappa}(\mathbb{R}) \times V_{p/2}(\Delta, \mathbb{R}^m) \times V_{p/2}(\Delta, \mathbb{R}^m)$ is a Gaussian measure with mean 0. Hence, the open ball centered at 0 has positive probability.

2. By Lemma 2.3,

$$E[\|C_{\mathbf{w}',w_{d+1}}\|_{p/2}^{p/2} \|w_{d+1}\|_{p,\kappa} < \varepsilon, \|C_{w_{d+1},\mathbf{z}}\|_{p/2} < a, \|C_{\mathbf{z},w_{d+1}}\|_{p/2} < a] \tag{3.21}$$
$$\leq C(p,\kappa) \cdot \varepsilon^{p/2},$$

where $C(p,\kappa)$ is a positive constant which depends only on $p$ and $\kappa$. By Chebyshev's inequality and Lemma 2.2.2, we obtain (3.16).

3. By (3.16),

$$1 - (\alpha_{a,a,\mathbf{z}} r)^2$$
$$\leq \mu_{\varepsilon,a,\mathbf{z}}(\|C_{\mathbf{w}',w_{d+1}}\|_{p/2} < a, \|C_{w_{d+1},\mathbf{w}'}\|_{p/2} < a)$$
$$(3.22) \quad = \int_{\mathbf{W}'} \left( \int_{W_{d+1}} \chi_{[0,a)}(\max\{\|C_{\mathbf{w}',w_{d+1}}\|_{p/2}, \right.$$
$$\left. \|C_{w_{d+1},\mathbf{w}'}\|_{p/2}\}) d\mu_{\varepsilon,a,\mathbf{z}}(w_{d+1}) \right) d\mu(\mathbf{w}')$$
$$= \int_{\mathbf{W}'} \mu_{\varepsilon,a,\mathbf{z}}(W_{a,(\mathbf{w}',\mathbf{z}),d+1}) d\mu(\mathbf{w}').$$

Hence, by Lemma 5.3 in [4], we get $\mu(V_{p,r,a,\mathbf{z},\varepsilon}) \geq 1 - r\alpha_{a,a,\mathbf{z}}$. Thus,

$$\mu(V_{p,r,a,\mathbf{z},\varepsilon} \cap \{\mathbf{w}' \in \mathbf{W}' | \|\overline{\mathbf{w}}'\|_{C^p} < a, \|C_{\mathbf{w}',\mathbf{z}}\|_{p/2} < a, \|C_{\mathbf{z},\mathbf{w}'}\|_{p/2} < a\})$$
$$(3.23) \quad \geq \mu(V_{p,r,a,\mathbf{z},\varepsilon}) - (1 - \alpha_{a,a,\mathbf{z}})$$
$$\geq (1-r)\alpha_{a,a,\mathbf{z}}.$$

4. Since $r < \frac{1}{3}$, $\mu_{\varepsilon,a,\mathbf{z}}(W_{a,(\mathbf{w}'_i,\mathbf{z}),d+1}) \geq \frac{2}{3}$. So it holds that

$$\mu_{\varepsilon,a,\mathbf{z}}(W_{a,(\mathbf{w}'_1,\mathbf{z}),d+1} \cap W_{a,(\mathbf{w}'_2,\mathbf{z}),d+1}) \geq \frac{1}{3},$$

which implies (3.20). □

The assertion that $W_{a,(\mathbf{w}',\mathbf{z}),d+1}$ is an $H$-convex set in $W_{d+1}$ and implies the following result which is a key to proving Lemma 2.5; refer to [6]. We denote $W_{a,(\mathbf{w}',\mathbf{z})}$, for simplicity, instead of $W_{a,(\mathbf{w}',\mathbf{z}),d+1}$.

LEMMA 3.4. *Let* $d\mu_{a,(\mathbf{w}',\mathbf{z})} = d\mu|_{W_{a,(\mathbf{w}',\mathbf{z})}}/(\mu(W_{a,(\mathbf{w}',\mathbf{z})}))$. *Let* $\mathcal{E}_{a,(\mathbf{w}',\mathbf{z})}$ *be the Dirichlet form on* $W_{a,(\mathbf{w}',\mathbf{z})}$. *Then for any* $f \in \mathrm{D}(\mathcal{E}_{W_{a,(\mathbf{w}',\mathbf{z})}})$ *the following LSI and PI hold:*

$$\int_{W_{a,(\mathbf{w}',\mathbf{z})}} f^2(w) \log(f^2(w)/\|f\|_{L^2(W_{a,(\mathbf{w}',\mathbf{z})})}^2) d\mu_{W_{a,(\mathbf{w}',\mathbf{z})}}(w)$$

(3.24)
$$\leq 2 \int_{W_{a,(\mathbf{w}',\mathbf{z})}} |Df(w)|_H^2 \, d\mu_{W_{a,(\mathbf{w}',\mathbf{z})}},$$

(3.25)
$$\iint_{W_{a,(\mathbf{w}',\mathbf{z})} \times W_{a,(\mathbf{w}',\mathbf{z})}} (f(w) - f(w'))^2 \, d\mu(w) \, d\mu(w')$$
$$\leq 2\mu(W_{a,(\mathbf{w}',\mathbf{z})}) \int_{W_{a,(\mathbf{w}',\mathbf{z})}} |Df(w)|_H^2 \, d\mu(w).$$

We prove Lemma 2.5 by using Lemma 3.1.

PROOF OF LEMMA 2.5. Recall the notation which we used in the proof of Lemma 2.6. First, we prove the case where $d=2$. Below we denote $U_{a,\mathbf{z}}$ by $U$ simply. Also, we use the notation in Lemma 3.1. In the present case, $Y_1 = W_1$, $Y_2 = W_2$, $m_i$ is the Wiener measure on $W_i$ and $m$ is the Wiener measure on $\mathbf{W}^2 = W_1 \times W_2$. In this case, $U_i = \{w_i \in W_i | \|C_{w_i,\mathbf{z}}\|_{p/2} < a, \|C_{\mathbf{z},w_i}\|_{p/2} < a, \|w_i\|_p < a\}$. Note that for $w_1 \in U_1$, $U_{w_1} = W_{a,(w_1,z),2}$. By Lemma 3.4, PI holds on $U_{w_1}$. Also PI holds on $U^{w_2}$ for the same reason. These statements imply that A1 holds for $U$. We prove A2. Let $V_{p,r,a,\mathbf{z},\varepsilon}$ be the set in (3.18). Then $m(U_1 \cap V_{p,r,a,\mathbf{z},\varepsilon}) \geq m(U_1)(1-r)$. Therefore, by (3.20), for sufficiently small $r$, $V_{p,r,a,\mathbf{z},\varepsilon}$ satisfies the property of A2. This completes the proof in the case of $d=2$. Now, we prove Lemma 2.5 in general dimension. We assume that Theorem 4.1 is valid in the case of $d$ dimension. We prove Theorem 4.1 in the case of $(d+1)$ dimension. We apply Lemma 3.1 in the case where $Y_1 = \mathbf{W}^d$, $Y_2 = W_{d+1}$ and $U = U_{a,\mathbf{z}}$. Note that the section of $U_{a,\mathbf{z}}$ by $\mathbf{w}'$ is nothing but $\mathbf{W}_{a,(\mathbf{w'z}),d+1} \subset W_{d+1}$. By Lemma 3.4, PI holds on the set. Also the section of $U_{a,\mathbf{z}}$ by $w_{d+1} \in V_{p,\kappa}(\mathbb{R})$ is $U_{a,(\mathbf{z},w_{d+1}),d}$. Therefore WPI holds on the set with a constant $\xi(\delta, a, \mathbf{z}, w_{d+1})$ by the assumption of induction. Clearly, we take this function to be measurable with respect to the variables $a$, $\delta$, $\mathbf{z}$ and $w_{d+1}$. These imply A1. By Lemma 3.3 parts 3 and 4, for any $\delta > 0$, we can find a subset $U_{a,\mathbf{z},d,\delta}$ of $U_{a,\mathbf{z},d}$ such that $\mu(U_{a,\mathbf{z},d} \setminus U_{a,\mathbf{z},d,\delta}) \leq \delta$ and, for any $\mathbf{w}', \mathbf{w}'' \in U_{a,\mathbf{z},d,\delta}$, $m_2(\mathbf{W}_{a,(\mathbf{w}',\mathbf{z})} \cap \mathbf{W}_{a,(\mathbf{w}'',\mathbf{z})}) > \beta(\delta, a, \mathbf{z})$, where $\beta(\delta, a, \mathbf{z})$ is a positive number. This implies A2. Consequently, we complete the proof. □

For $\mathbf{h} \in \mathbf{H}^d$ and $a > 0$, let

(3.26)
$$B_{a,\mathbf{h}} = \{\mathbf{w} \in \mathbf{X}^d | \|\overline{(\mathbf{w}-\mathbf{h})_2}\|_{p/2} < a,$$
$$\|C_{\mathbf{w}-\mathbf{h},\mathbf{h}}\|_{p/2} < a, \|C_{\mathbf{h},\mathbf{w}-\mathbf{h}}\|_{p/2} < a, \|\mathbf{w}-\mathbf{h}\|_p < a\}.$$

As a corollary of Lemma 2.5, we have the following lemma.

LEMMA 3.5. *For any $\mathbf{h} \in \mathbf{H}^d$ and $a > 0$, $B_{a,\mathbf{h}} = U_{a,\mathbf{h}} + \mathbf{h}$ a.s. and WPI holds on $B_{a,\mathbf{h}}$.*

PROOF. The equality $B_{a,\mathbf{h}} = U_{a,\mathbf{h}} + \mathbf{h}$ is obvious. Let $f \in \mathfrak{F}\mathfrak{C}_b^\infty$. Then applying WPI for $f(\mathbf{w}+\mathbf{h})$ on $U_{a,\mathbf{h}}$ and using the Cameron–Martin formula, we have

$$\int_{B_{a,\mathbf{h}}} (f(\mathbf{w}) - \mu(U_{a,\mathbf{h}})^{-1} E_\mu[f\rho_\mathbf{h}:B_{a,\mathbf{h}}])^2 \rho_\mathbf{h}(\mathbf{w})\mu(U_{a,\mathbf{h}})^{-1}\, d\mu(\mathbf{w})$$
(3.27)
$$\leq \xi(\delta,a,\mathbf{h}) \int_{B_{a,\mathbf{h}}} |Df(\mathbf{w})|_H^2 \mu(U_{a,\mathbf{h}})^{-1}\rho_\mathbf{h}(\mathbf{w})\, d\mu(\mathbf{w}) + \delta\|f\|_\infty^2,$$

where $\rho_\mathbf{h}(\mathbf{w}) = \exp((\mathbf{w},\mathbf{h}) - \|\mathbf{h}\|_{\mathbf{H}_d}^2/2)$. By Lemma 2.2 in [4], this completes the proof. $\square$

The following lemma shows that Lemma 2.4 is a stronger statement than Lemma 2.1.

LEMMA 3.6. 1. *For any* $\mathbf{h} \in \mathbf{H}^d$, *the following estimate holds:*

(3.28) $\|\overline{\mathbf{w}}_2 - \overline{\mathbf{h}}_2\|_{p/2} \leq \|\overline{(\mathbf{w}-\mathbf{h})_2}\|_{p/2} + 2\|C_{\mathbf{w}-\mathbf{h},\mathbf{h}}\|_{p/2} + \|\mathbf{h}\|_p \|\mathbf{w}-\mathbf{h}\|_p.$

2. *The following inclusion holds for any* $\varepsilon > 0$ *and* $\mathbf{h} \in \mathbf{H}^d$:

(3.29) $\quad B_{\varepsilon/(3+\|h\|_p),\mathbf{h}} \subset \{\mathbf{w} \in \mathbf{X}^d | \|\overline{\mathbf{w}}(\cdot,\cdot) - \overline{\mathbf{h}}(\cdot,\cdot)\|_{C^p} < \varepsilon\}.$

PROOF. Statement 1 follows from (2.26) immediately. Statement 2 follows from 1 immediately. $\square$

REMARK 3.7. At the moment, I do not know whether stronger PI or LSI hold on $U_{a,\mathbf{z}}$. Here, we prove that $U_{a,0}$ is not $H$ convex in the sense of [6] in the case of $d = 2$. This implies that the usual convexity criterion as in [6] is not applicable to $U_{a,0}$ at least. The proof is as follows. First note that the functional on $\mathbf{H}^2$ with values in $C^p(\Delta,\mathbb{R})$ such that $F(\mathbf{h})(s,t) = \int_s^t (h_1(u) - h_1(s))\dot{h}_2(u)\, du$, where $\mathbf{h} = (h_1, h_2)$, is not continuous in the topology of $V_p(\mathbb{R}^2)$; see [22]. Hence, there exists a sequence $\mathbf{h}_n = (h_1^n, h_2^n) \in \mathbf{H}^2$ such that $\limsup_{n\to\infty} \|\mathbf{h}_n\|_p < a$ and $\lim_{n\to\infty} \|F(\mathbf{h}_n)\|_{p/2} = \infty$. Set $\varphi_n = (h_1^n, 0) \in \mathbf{H}^2$, $\phi_n = (0, h_2^n) \in \mathbf{H}^2$ and $\eta_n = \phi_n - \varphi_n$. Then, for sufficiently large fixed $n$, there exists a small positive number $\varepsilon$ such that for almost all elements, $B_{\varepsilon,\varphi_n}, B_{\varepsilon,\varphi_n} + \eta_n \subset U_{a,0}$ by Lemma 2.2.1. However $B_{\varepsilon,\varphi_n} + \frac{1}{2}\eta_n \subset U_{a,0}^c$ for almost all elements. This shows that $U_{a,0}$ is not an $H$-convex set.

**4. Main theorem.** First we state our main theorem,

THEOREM 4.1. *Let $F$ be a real-valued function on $\mathbf{H}^d$ and assume that $F$ satisfies the following continuity condition. For any $R > 0$ and $\mathbf{h}_1, \mathbf{h}_2 \in \mathbf{H}^d$ with $\|\overline{\mathbf{h}}_1\|_{C^p}, \|\overline{\mathbf{h}}_2\|_{C^p} \leq R$, it holds that*

$$(4.1) \qquad |F(\mathbf{h}_1) - F(\mathbf{h}_2)| \leq C(R)\|\overline{\mathbf{h}}_1 - \overline{\mathbf{h}}_2\|_{C^p},$$

*where $C(R)$ is an increasing positive function of $R$. Then the following statements hold.*

1. *In the Hilbert space topology of $\mathbf{H}^d$, $F$ is a continuous function.*
2. *For any $\mathbf{w} \in \mathbf{X}^d$, $\lim_{n \to \infty} F(P_n \mathbf{w})$ converges. We denote the limit by $\widetilde{F}(\mathbf{w})$. Then $\widetilde{F}(\mathbf{w})$ is an $\mathbf{H}^d$-continuous function.*
3. *Let $\widetilde{U}_F = \{\mathbf{w} \in \mathbf{X}^d | \widetilde{F}(\mathbf{w}) > 0\}$ and $U_F = \{h \in \mathbf{H}^d | F(h) > 0\}$. Then $U_F \neq \varnothing$ is equivalent to $\mu(\widetilde{U}_F) > 0$. Also if $U_F$ is a connected set in $\mathbf{H}^d$, then WPI holds on $\widetilde{U}_F$.*

PROOF. By Lemma 2.2 and the assumption,

$$|F(\mathbf{h}_1) - F(\mathbf{h}_2)|$$
$$\leq C\left(C_1 \max_i \{(\|\mathbf{h}_i\|_{\mathbf{H}^d} + 1)\|\mathbf{h}_i\|_{\mathbf{H}^d}\}\right)(\|\mathbf{h}_1\|_{\mathbf{H}^d} + \|\mathbf{h}_2\|_{\mathbf{H}^d} + 1)\|\mathbf{h}_1 - \mathbf{h}_2\|_{\mathbf{H}^d}.$$

This proves statement 1. Since $\overline{P_n \mathbf{w}}$ converges in $G\Omega_p(\mathbb{R}^d)$, by the assumption of continuity, the convergence in part 2 is obvious. By (4.1), it holds that for any $\eta_1, \eta_2 \in \mathbf{X}^d$ with $\|\bar{\eta}_1\|_{C^p}, \|\bar{\eta}_2\|_{C^p} < R$,

$$|\widetilde{F}(\eta_1) - \widetilde{F}(\eta_2)| \leq C(R)\|\bar{\eta}_1 - \bar{\eta}_2\|_{C^p}.$$

By (2.26), we see the $\mathbf{H}^d$ continuity. Now we prove part 3. We see that the probability measure of $O_\varepsilon(\mathbf{h}) = \{\mathbf{w} \in \mathbf{X}^d | \|\overline{\mathbf{w}} - \overline{\mathbf{h}}\|_{C^p} \leq \varepsilon\}$ is positive by Lemmas 2.6 and 3.6 part 2. Assume that there exists $\mathbf{h} \in \mathbf{H}^d$ such that $F(\mathbf{h}) > 0$. Take $\mathbf{w} \in O_\varepsilon(\mathbf{h})$. Then by the assumption of $F$,

$$|\widetilde{F}(\mathbf{w}) - F(\mathbf{h})| \leq C(\|\overline{\mathbf{h}}\|_{C^p} + \varepsilon)\varepsilon.$$

This implies $\mu(\widetilde{F} > 0) > 0$. Conversely, we assume $\mu(\widetilde{F} > 0) > 0$. Then there exists $\mathbf{w} \in \mathbf{X}^d$ such that $\widetilde{F}(\mathbf{w}) > 0$. Then for sufficiently large $k$, $F(P_n \mathbf{w}) > 0$. Now we prove the latter half of statement 3. Take a countable dense set $\{\phi_i\}_{i=1}^\infty \subset U_F$ in the topology of $\mathbf{H}^d$. Let

$$B_{r,\phi_i,\mathbf{H}} = \{\mathbf{h} \in \mathbf{H}^d | \|\overline{(\mathbf{h} - \phi_i)_2}\|_{p/2} < r,$$
$$\|C_{\mathbf{h}-\phi_i,\phi_i}\|_{p/2} < r, \|C_{\phi_i,\mathbf{h}-\phi_i}\|_{p/2} < r, \|\mathbf{h} - \phi_i\|_p < r\}.$$

For each $\phi_i$, let $\{r_k^i\}_{k=1}^\infty$ be all positive rational numbers $r$ such that

$$(4.2) \qquad \inf\{F(\mathbf{h}) | \mathbf{h} \in B_{r,\phi_i,\mathbf{H}}\} > 0.$$

We use the following two claims; see (3.26) for the definition of $B_{r_k^i,\phi_i}$.

CLAIM 1. *We have*

(4.3) $$U_F = \bigcup_{i,k} B_{r_k^i, \phi_i, \mathbf{H}},$$

(4.4) $$\widetilde{U}_F = \bigcup_{i,k} B_{r_k^i, \phi_i} \quad a.s.$$

CLAIM 2. *The following two statements are equivalent:*

1. $B_{r_k^i, \phi_i, \mathbf{H}} \cap B_{r_k^j, \phi_j, \mathbf{H}} \neq \varnothing$;
2. $\mu(B_{r_k^i, \phi_i} \cap B_{r_k^j, \phi_j}) > 0$.

We can complete the proof of the theorem by these claims. By the connectivity assumption, we can change the order such that $\{B_{r_k^i, \phi_i, \mathbf{H}}\} = \{B_k\}_{k=1}^\infty$ and $(\bigcup_{k=1}^l B_k) \cap B_{l+1} \neq \varnothing$ for all $l \geq 1$. Let us denote by $\widetilde{B}_k$ the subset of $\widetilde{U}_F$ that corresponds to $B_k$. Then $\widetilde{U}_F = \bigcup_k \widetilde{B}_k$ and $\mu((\bigcup_{k=1}^l \widetilde{B}_k) \cap \widetilde{B}_{l+1}) > 0$. Therefore, by Theorem 6.10 in [4], WPI holds on $U_{\widetilde{F}}$. Now, we prove the claims. Noting the continuity of $F$, for sufficiently small $r$, we see that (4.2) holds, so the set on the right-hand side of (4.3) is a nonempty set. Take $\mathbf{h}_0 \in U_F$ and fix a rational number $0 < \varepsilon < 1$ such that

(4.5) $$\inf\{F(\mathbf{h}) | \mathbf{h} \in B_{\varepsilon, \mathbf{h}_0, \mathbf{H}}\} \geq \delta := \frac{F(\mathbf{h}_0)}{2}.$$

Then noting

(4.6) $$\begin{aligned}\overline{(\mathbf{h} - \phi_i)}(s,t)_2 &- \overline{(\mathbf{h} - \mathbf{h}_0)}(s,t)_2 \\ &= \overline{(\mathbf{h}_0 - \phi_i)}(s,t)_2 + \int_s^t [(\mathbf{h}_0 - \phi_i)(u) - (\mathbf{h}_0 - \phi_i)(s)] \otimes d(\mathbf{h} - \mathbf{h}_0)(u) \\ &\quad - \int_s^t d(\mathbf{h} - \mathbf{h}_0)(u) \otimes [(\mathbf{h}_0 - \phi_i)(u) - (\mathbf{h}_0 - \phi_i)(s)] \\ &\quad + [(\mathbf{h} - \mathbf{h}_0)(t) - (\mathbf{h} - \mathbf{h}_0)(s)] \otimes [(\mathbf{h}_0 - \phi_i)(t) - (\mathbf{h}_0 - \phi_i)(s)],\end{aligned}$$

by Lemma 2.2, we get

(4.7) $$\begin{aligned}\|\overline{(\mathbf{h} - \mathbf{h}_0)}_2\|_{p/2} &\leq \|\overline{(\mathbf{h} - \phi_i)}_2\|_{p/2} \\ &\quad + 6\|\mathbf{h}_0 - \phi_i\|_{\mathbf{H}^d}(\|\mathbf{h} - \phi_i\|_p + \|\phi_i - \mathbf{h}_0\|_{\mathbf{H}^d}),\end{aligned}$$

(4.8) $\|C_{\mathbf{h}-\mathbf{h}_0, \mathbf{h}_0}\|_{p/2} \leq \|C_{\mathbf{h}-\phi_i, \phi_i}\|_{p/2} + \|\mathbf{h}_0 - \phi_i\|_{\mathbf{H}^d}(\|\mathbf{h} - \phi_i\|_p + 2\|\mathbf{h}_0\|_p),$

(4.9) $\|\mathbf{h} - \mathbf{h}_0\|_p \leq \|\mathbf{h} - \phi_i\|_p + \|\phi_i - \mathbf{h}_0\|_{\mathbf{H}^d}.$

Hence, for $\phi_i$ and $\mathbf{h}$ with $\|\phi_i - \mathbf{h}_0\|_{\mathbf{H}^d} \leq \varepsilon_2$ and $\mathbf{h} \in B_{\varepsilon_1, \phi_i, \mathbf{H}}$, we have

(4.10) $$\|\overline{(\mathbf{h} - \mathbf{h}_0)}_2\|_{p/2} \leq \varepsilon_1 + 6\varepsilon_2(\varepsilon_1 + \varepsilon_2),$$

(4.11) $$\|C_{\mathbf{h}-\mathbf{h}_0,\mathbf{h}_0}\|_{p/2} \leq \varepsilon_1 + \varepsilon_2(\varepsilon_1 + 2\|\mathbf{h}_0\|_p),$$

(4.12) $$\|\mathbf{h} - \mathbf{h}_0\|_p \leq \varepsilon_1 + \varepsilon_2.$$

Hence, it holds that $\inf\{F(\mathbf{h})|\mathbf{h} \in B_{\varepsilon_1,\phi_i,\mathbf{H}}\} \geq \delta$ for $\phi_i$ above, $\varepsilon_1 = \frac{1}{6}\varepsilon$ and $\varepsilon_2 = \varepsilon/(36(1 + 2\|\mathbf{h}_0\|_p))$. By applying (4.7), (4.8) and (4.9) to the case where $\mathbf{h} = \phi_i$, we see that $\mathbf{h}_0 \in B_{\varepsilon_1,\phi_i,\mathbf{H}}$ for the same $\varepsilon_1$, $\varepsilon_2$ and $\phi_i$. This proves (4.3). Now we prove (4.4). Take $\eta \in \widetilde{U}_F$ and choose $\delta > 0$ and $R > 0$ such that $\widetilde{F}(\eta) > \delta$ and $\|\bar{\eta}\|_{C^p} \leq R < \infty$. Then by the definition of $\widetilde{F}$, there exists $\varepsilon > 0$ such that $\widetilde{F}(\mathbf{w}) \geq \delta/2$ holds for all $\mathbf{w}$ with $\|\overline{\mathbf{w}} - \bar{\eta}\|_{C^p} < \varepsilon$. Also there exists $N \in \mathbb{N}$ such that for all $n \geq N$, $\|\overline{P_n\eta} - \bar{\eta}\|_{C^p} < \varepsilon/2$. By Lemma 2.4, there exists $l > N$ such that $\eta \in B_{\varepsilon/(8(3+\|\eta\|_p)),P_l\eta}$. Choose $\phi_i$ such that $\|\phi_i - P_l\eta\|_{\mathbf{H}^d}$ is sufficiently small. Then, applying (4.6) to the case where $\mathbf{h} = \eta$ and $\mathbf{h}_0 = P_l\eta$, we have $\eta \in B_{\varepsilon/(4(3+\|\phi_i\|_p)),\phi_i}$. Also, by Lemma 3.6.2,

(4.13) $$B_{\varepsilon/(4(3+\|\phi_i\|_p)),\phi_i} \subset \{\mathbf{w} \in \mathbf{X}^d | \|\overline{\mathbf{w}} - \bar{\phi}_i\|_{C^p} < \varepsilon/4\}$$
$$\subset \{\mathbf{w} \in \mathbf{X}^d | \|\overline{\mathbf{w}} - \overline{P_l\eta}\|_{C^p} < \varepsilon/2\} \subset \widetilde{U}_F.$$

Assume (4.2). We need to prove that $B_{r,\phi_i} \subset \widetilde{U}_F$. Take $\mathbf{w} \in B_{r,\phi_i}$. Applying (2.26) to the case where $\varphi = \mathbf{w} - \phi_i$ and $\phi = P_n\mathbf{w} - \phi_i$, we have $P_n\mathbf{w} \in B_{r,\phi_i,\mathbf{H}}$ for all sufficiently large $n$. This implies $F(P_n\mathbf{w}) \geq \delta$ and $\widetilde{F}(\mathbf{w}) \geq \delta$. This proves (4.4). Now we prove Claim 2. Assume $B_{r_1,\phi_1,\mathbf{H}} \cap B_{r_2,\phi_2,\mathbf{H}} \neq \varnothing$. Then there exists $\mathbf{h} \in \mathbf{H}^d$ such that $\|\overline{\phi_i - \mathbf{h}}\|_{C^p} < r_i$, $\|C_{\phi_i-\mathbf{h},\phi_i}\|_{p/2} < r_i$ and $\|C_{\phi_i,\phi_i-\mathbf{h}}\|_{p/2} < r_i$. By (2.26),

(4.14)
$$C_{\phi_i-\mathbf{w},\phi_i-\mathbf{w}}(s,t) - C_{\phi_i-\mathbf{h},\phi_i-\mathbf{h}}(s,t)$$
$$= \overline{(\mathbf{h}-\mathbf{w})}(s,t)_2 + \int_s^t ((\mathbf{h}-\mathbf{w})(u) - (\mathbf{h}-\mathbf{w})(s)) \otimes d\phi_i(s)$$
$$- \int_s^t d\phi_i(u) \otimes ((\mathbf{h}-\mathbf{w})(u) - (\mathbf{h}-\mathbf{w})(s))$$
$$- \int_s^t ((\mathbf{h}-\mathbf{w})(u) - (\mathbf{h}-\mathbf{w})(s)) \otimes d\mathbf{h}(u)$$
$$+ \int_s^t d\mathbf{h}(u) \otimes ((\mathbf{h}-\mathbf{w})(u) - (\mathbf{h}-\mathbf{w})(s))$$
$$+ ((\phi_i - \mathbf{h})(t) - (\phi_i - \mathbf{h})(s)) \otimes ((\mathbf{h}-\mathbf{w})(t) - (\mathbf{h}-\mathbf{w})(s)).$$

Therefore, by Lemma 2.2, $\mathbf{w} \in B_{r_i,\phi_i}$ holds for $\mathbf{w}$ such that $\|\overline{\mathbf{h}-\mathbf{w}}\|_{C^p}$ is sufficiently small. This implies $\mu(B_{r_1,\phi_1} \cap B_{r_2,\phi_2}) > 0$. Next we assume $\mu(B_{r_1,\phi_1} \cap B_{r_2,\phi_2}) > 0$. Take $\mathbf{w} \in B_{r_1,\phi_1} \cap B_{r_2,\phi_2}$. Then by (4.14), $P_n\mathbf{w} \in B_{r_1,\phi_1,\mathbf{H}} \cap B_{r_2,\phi_2,\mathbf{H}}$. This completes the proof of Claim 2 and, hence, the theorem. $\square$

Finally we present examples.

EXAMPLE 1. Let $(M, g)$ be a $d$-dimensional compact Riemannian manifold. We consider an orthonormal frame bundle $\pi: O(M) \to M$. Let us take a metric connection and let $\{L_i\}_{i=1}^d$ be the corresponding canonical horizontal vector fields. Let us consider a Stratonovich stochastic differential equation

$$dr(t, u, \mathbf{w}) = \sum_{i=1}^d L_i(r(t, u, \mathbf{w})) \circ dw_i(t), \tag{4.15}$$

$$r(0) = u, \tag{4.16}$$

where $\pi(u) = x$ and $\mathbf{w} = (w_1, \ldots, w_d)$ denotes the $d$-dimensional Brownian motion. Let $X(t, x, \mathbf{w}) = \pi(r(t, \mathbf{w}))$, where $X(t, x, \mathbf{w})$ is a Brownian motion whose generator is $\Delta/2$, where $\Delta$ is the Laplace–Beltrami operator. Let us denote by $\eta(t, u, \mathbf{h})$ the solution to the ordinary differential equation which is given by replacing $\mathbf{w}$ with $\mathbf{h} \in \mathbf{H}^d$. We denote $\xi(t, x, \mathbf{h}) = \pi(\eta(t, u, \mathbf{h}))$. Let $y \in M$ and consider a small open geodesic ball $B_\varepsilon(y)$ centered at $y$ with radius $\varepsilon$. Let us consider a smooth function $\varphi$ on $M$ such that $\varphi(z) > 0$ holds if and only if $z \in B_\varepsilon(y)$. Then $F(\mathbf{h}) = \varphi(\xi(1, x, \mathbf{h}))$ satisfies the continuity assumption in Theorem 4.1 and $\widetilde{F}(\mathbf{w}) = \varphi(X(1, x, \mathbf{w}))$; see [17], Theorem 6.2.2, and [18], Proposition 6.2.2. Actually, the continuity result in Theorem 6.2.2 in [18] is a stronger statement than we require in our theorem. Note that the connectivity of $U_F$ is equivalent to the simply connectedness of $M$. Therefore, WPI holds on $U_{\widetilde{F}} = \{\mathbf{w} \in \mathbf{W}^d | X(1, x, \mathbf{w}) \in B_\varepsilon(y)\}$ if $M$ is simply connected.

EXAMPLE 2. Next suppose that $M$ is isometrically embedded into $\mathbb{R}^N$ and let $P(x): \mathbb{R}^N \to T_x M$ be the projection operators. Consider a gradient Brownian system

$$dX(t, x, \mathbf{w}) = P(X(t, x, \mathbf{w})) \circ d\mathbf{w}(t), \tag{4.17}$$

$$X(0, x, \mathbf{w}) = x. \tag{4.18}$$

In this case, $\mathbf{w} \in \mathbf{W}^N$. Suppose $M$ is simply connected. Then by the same argument as above, we have that WPI holds on $\{\mathbf{w} \in \mathbf{W}^N | X(1, x, \mathbf{w}) \in B_\varepsilon(y)\}$. By this, Lemmas 4.2–4.5 in [2] and Lemma 5.1 in [1], we see that WPI holds on the subset $\{\gamma \in P_x(M) | \gamma(1) \in B_\varepsilon(y)\}$ with natural Dirichlet form. Also, it is easy to prove that WPI holds for the natural Dirichlet forms on any open connected sets on $P_x(M)$ by a similar argument. Refer to [2] for the definition of Dirichlet forms.

EXAMPLE 3. We present an example for loop space. Suppose that $M$ is a compact Lie group $G$ with biinvariant Riemannian metric. Let $L_e(G) = C([0, 1] \to G | \gamma(0) = \gamma(1) = e)$, where $e$ denotes the unit element. By using the $H$ derivative $D_h F(\gamma) = \lim_{\varepsilon \to 0} \varepsilon^{-1}(F(e^{\varepsilon h(\cdot)} \gamma(\cdot)) - F(\gamma))$, we can define a Dirichlet form; see [9]. By using the tubular neighborhood, the retract map in [9], the conclusion in Example 1 above and Lemma 2.2 in [4], we see that WPI holds on $L_e(G)$ if $G$ is simply connected. We will study general cases for loop spaces over Riemannian manifolds in forthcoming papers.

**Acknowledgments.** I am grateful to Professor Terry Lyons for communicating his results and sending me a PDF file of his manuscript [18] before the publication. Professor Zhonming Qian and the referee pointed out a mistake in the proof of Lemma 2.3.1 in the first version and gave valuable comments on the definition of $U_{a,\mathbf{z}}$. I thank both of them.

## REFERENCES

[1] AIDA, S. (1998). Uniform positivity improving property, Sobolev inequality and spectral gaps. *J. Funct. Anal.* **158** 152–185. MR1641566

[2] AIDA, S. (1998). Differential calculus on path and loop spaces, II. Irreducibility of Dirichlet forms on loop spaces. *Bull. Sci. Math.* **122** 635–666. MR1668546

[3] AIDA, S. (2000). On the irreducibility of Dirichlet forms on domains in infinite dimensional spaces. *Osaka J. Math.* **37** 953–966. MR1809914

[4] AIDA, S. (2001). An estimate of the gap of spectrum of Schrödinger operators which generate hyperbounded semigroups. *J. Funct. Anal.* **185** 474–526. MR1856275

[5] BOGACHEV, V. I. (1998). *Gaussian Measures.* Amer. Math. Soc., Providence, RI. MR1642391

[6] FEYEL, D. and ÜSTÜNEL, A. S. (2000). The notion of convexity and concavity on Wiener space. *J. Funct. Anal.* **176** 400–428. MR1784421

[7] GONG, F., RÖCKNER, R. and WU, L. (2001). Poincaré inequality for weighted first order Sobolev spaces on loop spaces. *J. Funct. Anal.* **185** 527–563. MR1856276

[8] GONG, F. and WU, L. (2000). Spectral gap of positive operators and applications. *C. R. Acad. Sci. Paris Sér. I Math.* **331** 983–988. MR1809440

[9] GORSS, L. (1991). Logarithmic Sobolev inequalities on loop groups. *J. Funct. Anal.* **102** 268–313. MR1140628

[10] GROSS, L. (1993). Uniqueness of ground states for Schrödinger operators over loop groups. *J. Funct. Anal.* **112** 373–441. MR1213144

[11] HINO, M. (2000). Exponential decay of positivity preserving semigroups on $L^p$. *Osaka J. Math.* **37** 603–624. MR1789439

[12] ITÔ, K. and NISIO, M. (1968). On the convergence of sums of independent Banach space valued random variables. *Osaka J. Math.* **5** 35–84. MR235593

[13] KUSUOKA, S. (1991). Analysis on Wiener spaces, I. Nonlinear maps. *J. Funct. Anal.* **98** 122–168. MR1111196

[14] KUSUOKA, S. (1992). Analysis on Wiener spaces, II. Differential forms. *J. Funct. Anal.* **103** 229–274. MR1151548

[15] LEDOUX, M., LYONS, T. and QIAN, Z. (2002). Lévy area of Wiener processes in Banach spaces. *Ann. Probab.* **30** 546–578. MR1905851

[16] LEDOUX, M., QIAN, Z. and ZHANG, T. (2002). Large deviations and support theorem for diffusions via rough paths. *Stochastic Process. Appl.* **102** 265–283. MR1935127

[17] LYONS, T. (1998). Differential equations driven by rough signals. *Rev. Mat. Iberoamericana* **14** 215–310. MR1654527

[18] LYONS, T. and QIAN, Z. (2002). *System Control and Rough Paths*. Oxford Univ. Press. MR2036784
[19] MATHIEU, P. (1998). Quand l'inégalite log-Sobolev implique l'inégalite de trou spectral. *Séminaire de Probabilités XXXII. Lecture Notes in Math.* **1686** 30–35. Springer, Berlin. MR1651227
[20] MATHIEU, P. (1998). On the law of the hitting time of a small set by a Markov process. Preprint. MR1666883
[21] RÖCKNER, M. and WANG, F.-Y. (2001). Weak Poincaré inequalities and $L^2$-convergence rates of Markov semigroups. *J. Funct. Anal.* **185** 564–603. MR1856277
[22] SUGITA, H. (1989). Hu–Meyer's multiple Stratonovich integral and essential continuity of multiple Wiener integral. *Bull. Sci. Math.* **113** 463–474. MR1029620

DEPARTMENT OF MATHEMATICAL SCIENCE
GRADUATE SCHOOL OF ENGINEERING SCIENCE
OSAKA UNIVERSITY
TOYONAKA 560-8531
JAPAN
E-MAIL: aida@sigmath.es.osaka-u.ac.jp
URL: www.sigmath.es.osaka-u.ac.jp/˜aida/index.html